\documentclass[10pt]{article} 
\usepackage{amsmath}
\usepackage{amsfonts}
\begin{document}
\def\AA{A}%{\mbox{\boldmath $A$}}
\def\FF{F}%{\mbox{\boldmath $F$}}
\def\CC{C}%{\mbox{\boldmath $C$}}
\def\PP{P}%{\mbox{\boldmath $P$}}
\def\QQ{Q}%{\mbox{\boldmath $Q$}}
\def\RR{R}%{\mbox{\boldmath $R$}}
\def\ff{f}%{\mbox{\boldmath $f$}}
\def\aa{a}%{\mbox{\boldmath $a$}}
\def\gg{g}%{\mbox{\boldmath $g$}}
\def\xx{x}%{\mbox{\boldmath $x$}}
\def\DD{\mathcal{D}}
\def\non{ \nonumber }
\def\lar{\longrightarrow}
\def\llar{\longleftarrow}
\def\({\left(}
\def\){\right)}
%%%%%%%%%%%%%%%%%%%%%%%%%%%%%%%%%%%%%%%%%%%%%%%%%%%%%%%%%%%%%%%%%%
\newcommand{\mapup}[1]{\Big\uparrow
\rlap{$\vcenter{\hbox{$\scriptstyle#1$}}$}}
%%%%%%%%%%%%%%%%%%%%%%%%%%%%%%%%%%%%%%%%%%%%%%%%%%%%%%%%%%%%%%%%%%%%%
%\rightline{LPTHE-}
\vskip 1cm
\centerline{\LARGE Euler Characteristics of Theta Divisors }
\bigskip
\centerline{\LARGE  of Jacobians for Spectral Curves.}
\vskip 2cm
\centerline{\large 
A. Nakayashiki ${}^a$ and F.A. Smirnov ${}^b$
\footnote[0]{Membre du CNRS}}
\vskip1cm
\centerline{ ${}^a$ Faculty of
Mathematics, Kyushu University} 
\centerline{Ropponmatsu 4-2-1, Fukuoka 810-8560, Japan}
\bigskip
\centerline{ ${}^b$ Laboratoire de Physique Th\'eorique et Hautes
Energies \footnote[1]{\it Laboratoire associ\'e au CNRS.}}
\centerline{ Universit\'e Pierre et Marie Curie, Tour 16 1$^{er}$
		\'etage, 4 place Jussieu}
\centerline{75252 Paris cedex 05-France}
\vskip2cm
\noindent
{\bf Abstract.} 
We show how to calculate the Euler characteristic
of an affine Jacobi variety of a spectral curve
from its defining equations.

\newpage

\section{Introduction}

Consider an algebraic curve $X$, its Jacobian $J(X)$
and the theta divisor $\Theta\subset J(X)$.
The non-compact variety $J(X)-\Theta $ is known to be an affine algebraic 
variety, which we call the affine Jacobian variety of $X$.

In the paper \cite{nasm} we consider, following \cite{mum},
the algebraic model of an affine hyperelliptic Jacobian variety.
Under certain strong assumptions we calculated the topological Euler
characteristic and dimensions of singular cohomology groups of it.
More precisely we made a conjecture on the explicit forms of cohomologies
as the irreducible representations of the symplectic group.
Later one of the authors \cite{naka} was able to prove these conjectures. 
In this proof certain very particular properties of hyperelliptic
curves were used.

The algebraic model is available
for the affine Jacobians of more general algebraic curves of the form:
\begin{align}
w^N+t_1(z)w^{N-1}+\cdots +t_N(z)=0\ ,
\label{cur}
\end{align}
where $\text{deg}(t_j)\leq nj-1$ for some $n$ 
and the leading coefficient of $t_N$ does not vanish
(the corresponding construction is explained in detail in \cite{Beauv,ts}).
However, the cohomologies
of the affine Jacobians of these curves seem much more complicated
than in the hyperelliptic case.
Up to now we were not able to conjecture their explicit form. 
Thus even a partial information about those cohomologies is important.

In this paper we shall show that under a certain assumption 
(Assumption 1 in section 3)
the topological Euler characteristic $\chi (J(X)-\Theta)$
of the affine Jacobian of a curve $X$ from (\ref{cur}) can be calculated.
It then gives the Euler characteristic of the theta divisor $\Theta$
by the relation:
$$
\chi (\Theta)=-\chi (J(X)-\Theta).
$$

The paper is organized as follows. In the second section
we present a rather abstract scheme into which the algebraic
models of affine Jacobians can be put. We consider a certain family of 
affine algebraic sets parametrized by points of a complex affine space.
In the third section we consider the special case when the parameter
takes the origin. It is the only case when the affine ring becomes graded.
In this case the Euler characteristic can be easily 
determined in terms of the character of the affine ring. 
We define a double complex which permits to
identify the dimensions of cohomologies in question 
with the dimensions of the spaces of cochains in the minimal 
free resolution of the affine ring.
In the fourth section we show that the Euler characteristic of the variety 
at a generic point in the parameter space
coincides with the Euler characteristic of the variety at the origin
by considering certain deformation of the previous construction. 
This is the main result of this paper.
We apply this result to the family of the affine Jacobians of 
spectral curves (\ref{cur})
studied in \cite{ts} and get the formula for $\chi(\Theta)$.
In the last section we shall show that our formula for $\chi(\Theta)$
produces correct values in the case of genus $3$ and $4$. 
The appendix contains a proof of the proposition given in section two.

\section{General setting}

Consider the affine space $\mathbb{C}^m$ with coordinates
$(\aa _1,\cdots ,\aa _m)$. Consider $m-g$ polynomials
$\ff_j(\aa _1,\cdots ,\aa _m)$ and the affine algebraic set
$\mathcal{J}_{f^0}$ defined by the equations
\begin{align}
&\ff _1(\aa _1,\cdots ,\aa _m)=f_1^0,\non\\
&\qquad\cdots\non\\
&\ff_{n-g}(\aa _1,\cdots ,\aa _m)=f_{m-g}^0,
\non
\end{align}
where $f_{j}^0$ are complex numbers.
We assume that $\mathcal{J}_{0}=\{p\in {\mathbb C}^m| f_j(p)=0\,\, \forall j\}$
is $g$ dimensional, that is, the maximum of the dimensions of irreducible
components is $g$.
In the case of the family of the affine Jacobians of 
spectral curves (\ref{cur}) 
this assumption is verified by using the explicit parametrization of 
$\mathcal{J}_{0}$ \cite{ts,Beauv}.
We assign a positive degree $\text{deg}(\aa _j)$ to each 
$\aa _j$ and assume that
the polynomial $\ff _j$ is homogeneous of degree $\text{deg}(\ff _j)$ 
for all $j$.
We shall consider several rings:
\begin{eqnarray}
&&
\AA =\mathbb{C}[\aa _1,\cdots ,\aa _m],
\quad
\FF =\mathbb{C}[\FF _1,\cdots ,\FF_{m-g}],
\nonumber
\\
&&
\AA _0=\AA/\textstyle{\sum}(\ff _j\AA),
\quad
\AA _{f^0}=\AA/\textstyle{\sum}((\ff _j-f^0_j)\AA).
\nonumber
\end{eqnarray}
The rings $\FF$ and $\AA$ are polynomial rings of $m$ and $m-g$ variables
respectively.
The ring $\FF$ acts on $\AA$ by the multiplication of $F_j=f_j$.
We define $\text{deg}(\FF_j)=\text{deg}(\ff_j)$.
Then $\AA$, $\FF$ and $\AA_0$ become graded rings. 
But this is not the case for $\AA _{f^0}$ if some $f_j^0$ is not zero. 
The character of a graded vector space is
defined as the generating function of the dimensions of homogeneous
subspaces. It is obvious that the characters of $\AA$ and $\FF$ are given by
\begin{eqnarray}
&&
\text{ch}(\AA)=
\prod\limits _{j=1}^m\frac 1 {1-q^{\text{deg}(\aa _j)}},
\quad
\text{ch}(\FF) =\prod\limits _{j=1}^{m-g}\frac 1 {1-q^{\text{deg}(\ff _j)}}.
\nonumber
\end{eqnarray}

The determination of the character of $\AA_0$ is more involved.
But the result is simple.

\vskip 2mm\noindent
{\bf Proposition 1.} The character of $\AA_0$ is given by
\begin{eqnarray}
&&
\text{ch}(\AA _0)=\frac {\text{ch}(\AA)}{\text{ch}(\FF)}.
\label{charA0}
\end{eqnarray}
\vskip 2mm

The proof of this proposition is given in Appendix A.
For the case of the family of affine hyperelliptic Jacobians, 
this proposition were proved in \cite{nasm} by determining 
the ${\mathbb C}$-basis of $\AA_0$ explicitly. 
In a similar way to Proposition 2 in \cite{nasm} we have the
following corollary of this proposition.

\vskip 2mm\noindent
{\bf Corollary 1.}  As a graded $\FF$ module, $\AA$ is a free module:
$\AA \simeq \FF\otimes_{\mathbb{C}} \AA _0$.
\vskip 2mm

\noindent 
An important feature, which is special for the case related with the family 
of affine Jacobians, is that on $\AA$ one can define the action of $g$
commuting vector fields $D_1,\cdots ,D_g$. 
Thus we assume that there exist algebraically independent $g$ 
commuting vector fields $D_1,\cdots ,D_g$ acting on $\AA$ such that 
\begin{eqnarray}
&&
 D_j\ff _k=0\quad\forall j,k.
\label{vectorfield}
\end{eqnarray}
Moreover, we assume that, to every $D_j$ a positive degree is 
prescribed such that
$$
\text{deg}(D_j\xx )=\text{deg}(D_j)+\text{deg}(\xx ) 
$$
for every homogeneous $x\in \AA$.
Then the ring 
$\DD =\mathbb{C}[D_1,\cdots,D_g]$
becomes a graded ring and its character is given by
$$
\text{ch}(\DD )=\prod\limits _{j=1}^g\frac 1 {1-q^{\text{deg}(D _j)}}\ .
$$
By the condition (\ref{vectorfield}) 
the action of $\DD$ on $\AA$ descends to $\AA _0$ and $\AA _{f^0}$.

Let us introduce differentials $dz _j$ dual to the vector fields $D_j$ .
We set $V=\oplus_{j=1}^g\mathbb{C}dz_j$. We define the spaces of cochains by
$$ 
\CC ^k=\AA \otimes \wedge^k V
,\quad \CC ^k_0=\AA_0\otimes \wedge^k V
,\quad \CC ^k_{f^0}=\AA _{f^0}\otimes \wedge^k V\ .
$$
The differential $d=\sum D_j\otimes dz_j$ defines complexes
$(\CC^{\cdot},d)$, $(\CC^{\cdot}_0,d)$, $(\CC^{\cdot}_{f^0},d)$
which we shall study. 
We set $\text{deg}(dz _j)=-\text{deg}(D_j)$ in order that 
$\text{deg}(d)=0$. Then $\CC^k$ and $\CC^k_0$ become graded vector spaces.

\section{Cohomologies at the origin}
Consider the complex:
\begin{align}
&0\lar
{\CC ^{0}_0}
\stackrel{d}{\lar}
\ {\CC ^{1}_0}
\stackrel{d}{\lar }
\ \cdots
\stackrel{d}{\lar }
\ {\CC ^{g-1}_0}
\stackrel{d}{\lar }
\ {\CC ^{g}_0}
\stackrel{d}{\lar }\ 0.
\end{align}
This complex is graded. One easily finds its $q$-Euler characteristic
\cite{nasm}:
$$
\chi _q (\CC ^{\cdot}_0)=
\sum_{k=0}^g(-1)^k\text{ch}(\CC^k_0)=
(-1)^gq^{-\sum \text{deg}(D_j)}
\frac {\text{ch}(\AA _0)}{\text{ch}(\DD )}\ .
$$
Taking the limit $q\to 1$  we find
the Euler characteristic as a number:
\begin{align}
&\chi (\CC ^{\cdot}_0)=
(-1)^g\frac {\prod\limits _{j=1}^{m-g}\text{deg}(\ff _j)
\prod\limits _{j=1}^g\text{deg}(D _j)}
{\prod\limits _{j=1}^m\text{deg}(\aa _j)}.
\label{euler}
\end{align}
 
The fact that the Euler characteristic is finite does not
mean that the cohomologies are finite dimensional. 
However we shall adopt the following assumption:
\vskip 2mm
\noindent
{\bf Assumption 1.} $\hbox{dim} H^g(\CC^{\cdot}_0)<\infty$.
\vskip 2mm
\noindent
Unfortunately,
we do not know how to prove this assumption even in the case of the family of
affine Jacobians of (\ref{cur}).

The Assumption 1 together with simple grading arguments guarantee
that $\AA_0$ is a finitely generated $\DD$-module. 
Then it has the minimal
$\DD$ free resolution ({\it cf.} \cite{Eis} \S 19) of the form
\begin{eqnarray}
&&
0
\lar 
\DD^{b_{0}}
\stackrel{\varphi_{0}^{0}(D)}{\lar} 
\cdots
\lar
\DD^{b_{g-2}}
\stackrel{\varphi_{0}^{g-2}(D)}{\lar} 
\DD^{b_{g-1}}
\stackrel{\varphi_{0}^{g-1}(D)}{\lar} 
\DD^{b_g}
\lar 
\AA_0
\lar
0.
\label{resol1}
\end{eqnarray}
Here minimality condition is $\varphi_{0}^k(0)=0$ for all $k$.
Consider another complex, the Koszul cpmplex of $D_1$, ..., $D_g$:
\begin{eqnarray}
&&
0
\lar 
\DD
\stackrel{d}{\lar}
\cdots
\lar
\DD\otimes\wedge^{g-2}V
\stackrel{d}{\lar} 
\DD\otimes\wedge^{g-1}V
\stackrel{d}{\lar} 
\DD\otimes\wedge^{g}V
\lar 
\mathbb{C}
\lar
0,
\label{Koszul}
\end{eqnarray}
where $d$ is given by the same formula as before $d=\sum D_j\otimes dz_j$.
This complex is exact. For any non-negative integer $b$, tensoring $\DD^b$ 
(over $\DD$) to (\ref{Koszul}), we get the exact sequence
\begin{eqnarray}
&&
0
\lar 
\DD^b
\stackrel{d}{\lar } 
\cdots
\lar
\DD^b\otimes\wedge^{g-2}V
\stackrel{d}{\lar } 
\DD^b\otimes\wedge^{g-1}V
\stackrel{d}{\lar } 
\DD^b\otimes\wedge^{g}V
\lar 
\mathbb{C}^b
\lar
0.
\nonumber
\end{eqnarray}
Now we can consider the double complex
$$
\begin{array}{ccccccccccccc}
{}&{}&0&{}&0&{}&0&{}&0&{}&0&{}&{}
\\
{}&{}&\Big\uparrow&{}&\Big\uparrow&{}&\Big\uparrow&{}&\Big\uparrow&{}&
\Big\uparrow&{}&{}
\\
0&\lar&{\AA_0}&\stackrel{d}{\lar}&{\CC ^{1}_0}&\stackrel{d}{\lar}
\cdots\stackrel{d}{\lar}&{\CC ^{g-1}_0}&\stackrel{d}{\lar}&\CC ^{g}_0
&\lar&0&\lar&0
\\
{}&{}&\Big\uparrow&{}&\Big\uparrow&{}&\Big\uparrow&{}&\Big\uparrow
&{}&\Big\uparrow&{}&{}
\\
0&\lar&\DD^{b_g}&\stackrel{d}{\lar}&\DD^{b_g}\otimes V&
\stackrel{d}{\lar}\cdots\stackrel{d}{\lar}&
\DD^{b_g}\otimes \wedge^{g-1} V&\stackrel{d}{\lar}&\DD^{b_g}\otimes \wedge^g V
&\stackrel{d}{\lar}&{\mathbb{C}^{b_g}}&\lar&0
\\
{}&{}&\mapup{\varphi_{0}^{g-1}(D)\otimes 1}&{}&
\mapup{\varphi_{0}^{g-1}(D)\otimes 1}&{}&
\mapup{\varphi_{0}^{g-1}(D)\otimes 1}&{}&
\mapup{\varphi_{0}^{g-1}(D)\otimes 1}&{}&
\Big\uparrow&{}&{}
\\
0&\lar&\DD^{b_{g-1}}&\stackrel{d}{\lar}&\DD^{b_{g-1}}\otimes V&
\stackrel{d}{\lar}\cdots\stackrel{d}{\lar}&
\DD^{b_{g-1}}\otimes \wedge^{g-1} V&\stackrel{d}{\lar}&
\DD^{b_{g-1}}\otimes \wedge^g V&\stackrel{d}{\lar}&
{\mathbb{C}^{b_{g-1}}}&\lar&0
\\
{}&{}&\Big\uparrow&{}&\Big\uparrow&{}&\Big\uparrow&{}&\Big\uparrow&{}&
\Big\uparrow&{}&{}
\\
{}&{}&\vdots&{}&\vdots&{}&\vdots&{}&\vdots&{}&\vdots&{}&{}
\\
{}&{}&\mapup{\varphi_{0}^{0}(D)\otimes 1}&{}&\mapup{\varphi_{0}^{0}(D)\otimes 
1}&{}&\mapup{\varphi_{0}^{0}(D)\otimes 1}&{}&\mapup{\varphi_{0}^{0}(D)\otimes 1}
&{}&\Big\uparrow&{}&{}
\\
0&\lar&\DD^{b_0}&\stackrel{d}{\lar}&\DD^{b_0}\otimes V&
\stackrel{d}{\lar}\cdots\stackrel{d}{\lar}&
\DD^{b_0}\otimes \wedge^{g-1} V&\stackrel{d}{\lar}&
\DD^{b_0}\otimes \wedge^g V&\stackrel{d}{\lar}&{\mathbb{C}^{b_0}}&\lar&0
\\
{}&{}&\Big\uparrow&{}&\Big\uparrow&{}&\Big\uparrow&{}&\Big\uparrow&{}&
\Big\uparrow&{}&{}
\\
{}&{}&0&{}&0&{}&0&{}&0&{}&0{}&{}.
\\
\end{array}
$$ 
In this double complex each column except the last one is exact, 
each row except the first one is exact. 
The minimality condition means that maps of the last column
${\mathbb C}^{b_k}\lar{\mathbb C}^{b_{k+1}}$, $k\leq g-1$ are all zero maps.
Using this diagram one can easily prove that
\begin{eqnarray}
&&
H^k(\CC_0^{\cdot})\simeq {\mathbb C}^{b_k}.
\nonumber
\end{eqnarray}
\vskip5mm

\section{Cohomologies at a generic point}
We set $\DD_F=\FF\otimes_{{\mathbb C}}\DD$. 
Under Assumption 1, using Proposition 1,
we can prove that $\AA$ is a finitely generated $\DD_F$ module
in a similar manner to Proposition 4 in \cite{nasm}.
By the uniqueness of the minimal resolution 
({\it cf.} \cite{Eis}, the graded version of Theorem 20.2), 
(\ref{resol1}) implies the existence of the
minimal $\DD_F$-free resolution of $\AA$ of the form
\begin{eqnarray}
&&
0
\lar 
\DD_F^{b_{0}}
\stackrel{{\varphi}^{0}(D)}{\lar} 
\cdots
\lar
\DD_F^{b_{g-2}}
\stackrel{{\varphi}^{g-2}(D)}{\lar} 
\DD_F^{b_{g-1}}
\stackrel{{\varphi}^{g-1}(D)}{\lar} 
\DD_F^{b_g}
\lar 
\AA
\lar
0,
\nonumber
\end{eqnarray}
such that 
$
{\varphi}^{k}(D)\vert_{\ff _j=0}=\varphi_{0}^{k}(D).
$
The minimality condition is that 
$
{\varphi}^{k}(D)\vert_{\ff _j=0,D_j=0}=0.
$
Again consider the double complex:
$$
\begin{array}{ccccccccccccc}
{}&{}&0&{}&0&{}&0&{}&0&{}&0&{}&{}
\\
{}&{}&\Big\uparrow&{}&\Big\uparrow&{}&\Big\uparrow&{}&\Big\uparrow&{}&
\Big\uparrow&{}&{}
\\
0&\lar&{\AA}&\stackrel{d}{\lar}&{\CC ^{1}}&\stackrel{d}{\lar}
\cdots\stackrel{d}{\lar}&{\CC ^{g-1}}&\stackrel{d}{\lar}&\CC ^{g}
&\lar&0&\lar&0
\\
{}&{}&\Big\uparrow&{}&\Big\uparrow&{}&\Big\uparrow&{}&\Big\uparrow
&{}&\Big\uparrow&{}&{}
\\
0&\lar&\DD_{F}^{b_g}&\stackrel{d}{\lar}&\DD_{F}^{b_g}\otimes V&
\stackrel{d}{\lar}\cdots\stackrel{d}{\lar}&
\DD_{F}^{b_g}\otimes \wedge^{g-1} V&
\stackrel{d}{\lar}&
\DD_{F}^{b_g}\otimes \wedge^g V
&\lar&{\FF^{b_g}}&\lar&0
\\
{}&{}&\mapup{\varphi^{g-1}(D)\otimes 1}&{}&
\mapup{\varphi^{g-1}(D)\otimes 1}&{}&
\mapup{\varphi^{g-1}(D)\otimes 1}&{}&
\mapup{\varphi^{g-1}(D)\otimes 1}&{}&
\mapup{\varphi^{g-1}(0)\otimes 1}&{}&{}
\\
0&\lar&\DD_{F}^{b_{g-1}}&\stackrel{d}{\lar}&
\DD_{F}^{b_{g-1}}\otimes V&
\stackrel{d}{\lar}\cdots\stackrel{d}{\lar}&
\DD_{F}^{b_{g-1}}\otimes \wedge^{g-1} V&\stackrel{d}{\lar}&
\DD_{F}^{b_{g-1}}\otimes \wedge^g V
&\lar&{\FF^{b_{g-1}}}&\lar&0
\\
{}&{}&\Big\uparrow&{}&\Big\uparrow&{}&\Big\uparrow&{}&\Big\uparrow&{}&
\Big\uparrow&{}&{}
\\
{}&{}&\vdots&{}&\vdots&{}&\vdots&{}&\vdots&{}&\vdots&{}&{}
\\
{}&{}&
\mapup{\varphi^{0}(D)\otimes 1}&{}&
\mapup{\varphi_{0}^{0}(D)\otimes 1}&{}&
\mapup{\varphi^{0}(D)\otimes 1}&{}&
\mapup{\varphi_{0}^{0}(D)\otimes 1}
&{}&
\mapup{\varphi^{0}(0)\otimes 1}&{}&{}
\\
0&\lar&\DD_{F}^{b_0}&\stackrel{d}{\lar}&
\DD_{F}^{b_0}\otimes V&
\stackrel{d}{\lar}\cdots\stackrel{d}{\lar}&
\DD_{F}^{b_0}\otimes \wedge^{g-1} V&\stackrel{d}{\lar}&
\DD_{F}^{b_0}\otimes \wedge^g V
&\lar&{\FF^{b_0}}&\lar&0
\\
{}&{}&\Big\uparrow&{}&\Big\uparrow&{}&\Big\uparrow&{}&\Big\uparrow&{}&
\Big\uparrow&{}&{}
\\
{}&{}&0&{}&0&{}&0&{}&0&{}&0{}&{}.
\\
\end{array}
$$ 
In the diagram each column except the last one and each row except 
the first one are exact.
Since every map in the diagram is $\FF$-linear and every object in the diagram
is an $\FF$-free module, by tensoring 
${\mathbb C}=\FF/\sum (f_j-f_j^0)\FF$ over $\FF$, we get the double complex:
$$
\begin{array}{ccccccccccccc}
{}&{}&0&{}&0&{}&0&{}&0&{}&0&{}&{}
\\
{}&{}&\Big\uparrow&{}&\Big\uparrow&{}&\Big\uparrow&{}&\Big\uparrow&{}&
\Big\uparrow&{}&{}
\\
0&\lar&{\AA_{f^0}}&\stackrel{d}{\lar}&{\CC ^{1}_{f^0}}&\stackrel{d}{\lar}
\cdots\stackrel{d}{\lar}&{\CC ^{g-1}_{f^0}}&\stackrel{d}{\lar}&\CC ^{g}_{f^0}
&\lar&0&\lar&0
\\
{}&{}&\Big\uparrow&{}&\Big\uparrow&{}&\Big\uparrow&{}&\Big\uparrow
&{}&\Big\uparrow&{}&{}
\\
0&\lar&\DD^{b_g}&\stackrel{d}{\lar}&\DD^{b_g}\otimes V&
\stackrel{d}{\lar}\cdots\stackrel{d}{\lar}&
\DD^{b_g}\otimes \wedge^{g-1} V&\stackrel{d}{\lar}&
\DD^{b_g}\otimes \wedge^g V
&\lar&{\mathbb{C}^{b_g}}&\lar&0
\\
{}&{}&\mapup{\varphi_{f^0}^{g-1}(D)\otimes 1}&{}&
\mapup{\varphi_{f^0}^{g-1}(D)\otimes 1}&{}&
\mapup{\varphi_{f^0}^{g-1}(D)\otimes 1}&{}&
\mapup{\varphi_{f^0}^{g-1}(D)\otimes 1}&{}&
\mapup{\varphi_{f^0}^{g-1}(0)\otimes 1}&{}&{}
\\
0&\lar&\DD^{b_{g-1}}&\stackrel{d}{\lar}&
\DD^{b_{g-1}}\otimes V&
\stackrel{d}{\lar}\cdots\stackrel{d}{\lar}&
\DD^{b_{g-1}}\otimes \wedge^{g-1} V&\stackrel{d}{\lar}&
\DD^{b_{g-1}}\otimes \wedge^g V
&\lar&{\mathbb{C}^{b_{g-1}}}&\lar&0
\\
{}&{}&\Big\uparrow&{}&\Big\uparrow&{}&\Big\uparrow&{}&\Big\uparrow&{}&
\Big\uparrow&{}&{}
\\
{}&{}&\vdots&{}&\vdots&{}&\vdots&{}&\vdots&{}&\vdots&{}&{}
\\
{}&{}&\mapup{\varphi_{f^0}^{0}(D)\otimes 1}&{}
&\mapup{\varphi_{f^0}^{0}(D)\otimes 1}&{}
&\mapup{\varphi_{f^0}^{0}(D)\otimes 1}&{}
&\mapup{\varphi_{f^0}^{0}(D)\otimes 1}
&{}&
\mapup{\varphi_{f^0}^{0}(0)\otimes 1}&{}&{}
\\
0&\lar&\DD^{b_0}&\stackrel{d}{\lar}&
\DD^{b_0}\otimes V&
\stackrel{d}{\lar}\cdots\stackrel{d}{\lar}&
\DD^{b_0}\otimes \wedge^{g-1} V&\stackrel{d}{\lar}&
\DD^{b_0}\otimes \wedge^g V
&\lar&{\mathbb{C}^{b_0}}&\lar&0
\\
{}&{}&\Big\uparrow&{}&\Big\uparrow&{}&\Big\uparrow&{}&\Big\uparrow&{}&
\Big\uparrow&{}&{}
\\
{}&{}&0&{}&0&{}&0&{}&0&{}&0{}&{},
\\
\end{array}
$$ 
where
$$
\varphi_{f^0}^k(D)={\varphi}^k(D)\vert_{f_j=f_j^0}.
$$
In the diagram each column except the last one and each row except 
the first one are exact.
One can easily prove the following isomorphism from this diagram:
\begin{eqnarray}
&&
H^k(\CC_{f^0}^{\cdot})\simeq 
\hbox{Ker}(\varphi_{f^0}^k(0)\otimes 1)/
\hbox{Im}(\varphi_{f^0}^{k-1}(0)\otimes 1),
\quad
1\leq k\leq g,
\nonumber
\end{eqnarray}
where we set $\varphi_{f^0}^g(0)=0$.
It follows that
\begin{eqnarray}
\chi(\CC^{\cdot}_{f^0})&=&
\sum_{k=0}^g(-1)^k
\Big(
\text{dim}\big(\hbox{Ker}(\varphi_{f^0}^k(0)\otimes 1)\big)
-\text{dim}\big(\hbox{Im}(\varphi_{f^0}^{k-1}(0)\otimes 1)\big)
\Big)
\nonumber
\\
&=&
\sum_{k=0}^g(-1)^k
\Big(
\text{dim}\big(\hbox{Ker}(\varphi_{f^0}^k(0)\otimes 1)\big)
+\text{dim}\big(\hbox{Im}(\varphi_{f^0}^{k}(0)\otimes 1)\big)
\Big)
\nonumber
\\
&=&
\sum_{k=0}^g(-1)^kb_k=\chi(\CC^{\cdot}_0).
\nonumber
\end{eqnarray}

Thus the cohomologies at generic point can differ from 
the cohomologies at the origin. But the Euler characteristic
is the same. 
Summarizing the results we have
\vskip2mm

\noindent
{\bf Theorem.} For any $f^0$, $\chi(\CC^{\cdot}_{f^0})=\chi(\CC^{\cdot}_{0})$ 
and it is given by (\ref{euler}).
\vskip2mm

\noindent
Now we apply this theorem to the case studied in \cite{ts}.
Then for the curve (\ref{cur}), which we assume non-singular,
of genus $g=\frac 1 2 (N-1)(Nn-2)$ using the results of \cite{ts}
one finds:
\begin{align}
\chi (\Theta )=(-1)^{g-1}N^{N^2n-2N+1}
(Nn-1)^{N-1}(\Gamma (N))^2
%\prod\limits _{j=1}^N\frac {\Gamma (j)}{\Gamma (nN+j-1)}
\prod\limits _{j=1}^{N-1}\frac {\Gamma (j)}
{\Gamma (nN+j)}
\left(\frac {\Gamma\left (j\frac{nN-1}{N}\right)}
{\Gamma \left(\frac{j}{N}\right)}\right) ^2\ .
\label{main}
\end{align}
This formula is quite different from what we have for a generic Abelian
variety, $ \chi (\Theta )=(-1)^{g-1}g!$
(see \cite{nasm}).

The formula (\ref{main}) is still a conjecture since Assumption 1
is not proved in this special case.

\section{Examples}
For the hyperelliptic case $N=2$, the formula (\ref{main}) was
proved in \cite{naka}.
In this section we shall give other examples for which (\ref{main})
is verified.
\vskip2mm

\noindent
Example 1. $g=3$ case: Consider the case $N=4$, $n=1$ and $X$ is given by
$w^4+t_4(z)=0$, where $t_4(z)$ is a polynomial of degree $3$.
We assume that $t_4(z)$ does not have multiple zeros.
In this case genus $g$ of $X$ is $3$.
The curve $X$ is not hyperelliptic since $w$ has third order poles at infinity
and is holomorphic at all other points.
Thus $\Theta$ is non-singular and $\chi(\Theta)=6$.
Our formula (\ref{main}) recovers this.
\vskip3mm

\par
\noindent
Example 2. $g=4$ case:
Consider the case when $N=3$, $n=2$ and the curve $X$ is given by
$w^3+t_3(z)=0$. 
The polynomial $t_3(z)$ is of degree $5$,
we suppose that it does not have multiple zeros.
The genus $g$ of $X$ is $4$.
In this case the singularity of the theta divisor $\Theta$ consists
of one point, say $e$ ({\it cf.} \cite{FK}).
We have the Abel-Jacobi map $\pi$ from the symmetric products of $X$ 
to $J(X)$. By Riemann's theorem the restriction of $\pi$ to $S^3(X)$ gives
the surjective map
$$
\pi: S^3(X) \lar \Theta.
$$
Moreover the isomorphism
$$
S^3(X)-\pi^{-1}(e)\simeq \Theta-\{e\}
$$
holds. Since $\pi^{-1}(e)\simeq {\mathbb P}^1$, the Euler characteristic 
of $\Theta$ is given by
\begin{eqnarray}
&&
\chi(\Theta)=\chi(S^3(X))-\chi(\pi^{-1}(e))+1=-21.
\nonumber
\end{eqnarray}
Our formula (\ref{main}) precisely recovers this.
Notice that the curve $X$ is not hyperelliptic and 
is not generic among genus $4$ non-hyperelliptic curves.
\vskip5mm

Those examples support the validity of our assumptions.

For Jacobians of curves satisfying certain genericity conditions,
the formula of $\chi(\Theta)$ is given in \cite{PP}.
For genus 4 case, their formula produces $\chi(\Theta)=-22$.
This corresponds to curves whose $\Theta$ has exactly two singular points.
\vskip3mm

\noindent
{\bf \large Acknowledgements} \par
\noindent
We would like to thank Fu Baohua for comments on the appendix.
One of the authors (FS) was partly supported by INTAS
grant INTAS OPEN 97-01312.

\section{Appendix A}
We shall give a proof of Proposition 1.

It is obvious that (\ref{charA0}) holds if and only if 
$(f_1,\cdots,f_{m-g})$ is a regular sequence.
We set $I_j=\AA f_1+\cdots+\AA f_j$ and, for a subset $T\subset \AA$,
$Z(T)=\{p\in {\mathbb C}^m| f(p)=0 \quad\forall f\in T\}$.

We first notice that $f_{j+1}\notin I_j$ for any $j$. Because, otherwise
$\dim (\mathcal{J}_{0})>g$.
Suppose that $(f_1,\cdots,f_{m-g})$ is not a regular sequence.
Then there exist $i$ and $g\in \AA$ such that
\begin{eqnarray}
&&
f_{i+1}g\in I_i,
\quad
g\notin I_i.
\nonumber
\end{eqnarray}
Write $I_i$ as an intersection of homogeneous primary ideals ${\cal Q}_j$
in an irredundant way ({\it cf}. \cite{zs}, Ch IV, \S 4)
\begin{eqnarray}
&&
I_i={\cal Q}_1\cap\cdots\cap {\cal Q}_r.
\nonumber
\end{eqnarray}
This is possible since $f_j$'s are all homogeneous.
Then $f_{i+1}\notin \sqrt{{\cal Q}_j}$ for any $j$.
In fact suppose that $f_{i+1}\in \sqrt{{\cal Q}_{j_0}}$ for some $j_0$.
Then 
\begin{eqnarray}
&&
Z(f_{i+1})\supset Z(\sqrt{{\cal Q}_{j_0}})=Z({\cal Q}_{j_0}).
\label{inclusion}
\end{eqnarray}
Notice that 
\begin{eqnarray}
&&
Z(I_i)=Z({\cal Q}_1)\cup \cdots \cup Z({\cal Q}_r)
\nonumber
\end{eqnarray}
is the irreducible decomposition of $Z(I_i)$.
Since $Z(I_i)$ is the zero set of $i$ polynomials, 
$\dim (Z({\cal Q}_j))\geq m-i$ for any $j$.
Then the relation (\ref{inclusion}) implies that 
$\dim (Z(I_{i+1}))\geq \dim (Z({\cal Q}_{j_0}))\geq m-i$.
Since ${\cal Q}_j$'s and $f_j$'s are homogeneous,
the intersection 
$Z({\cal Q}_{j_0})\cap Z(f_{i+2})\cap\cdots\cap Z(f_{m-g})$
is non-empty.
Thus $\dim (Z(I_{m-g}))=\dim (\mathcal{J}_{0})>g$,
which contradicts the assumption.

On the other hand $g\notin {\cal Q}_{j_0}$ for some $j_0$, since
$g\notin I_i$. This means that $f_{i+1}g \in {\cal Q}_{j_0}$ and
$g\notin {\cal Q}_{j_0}$.
It follows that $f_{i+1}^k\in {\cal Q}_{j_0}$ for some $k\geq 2$, which
is impossible.

\end{document}